\documentclass{article}
\usepackage{graphicx}
\usepackage{amsmath}
\usepackage{epsfig}

\newtheorem{theorem}{Theorem}[section]

\newtheorem{lemma}[theorem]{Lemma}

\newenvironment{proof}[1][Proof]{\textbf{#1.} }{\ \rule{0.5em}{0.5em}}

\begin{document}

\title{Cubic Maximal Nontraceable Graphs\thanks{%
This material is based upon work supported by the National Research
Foundation under Grant number 2053752.}}
\author{Marietjie Frick, Joy Singleton \\
University of South Africa,\\
P.O. Box 392, Unisa, 0003,\\
South Africa.\\
\textbf{e-mail:} singlje@unisa.ac.za frickm@unisa.ac.za}
\date{ }
\maketitle

\begin{abstract}
We determine a lower bound for the number of edges of a 2-connected maximal
nontraceable graph, and present a construction of an infinite family of
maximal nontraceable graphs that realize this bound.


\end{abstract}

\section{\ Introduction}

We consider only simple, finite graphs $G$ and denote the vertex set and
edge set of $G$ by $V(G)$ and $E(G)$, respectively. The \textit{open neighbourhood} 
of a vertex $v$ in $G$ is the set 
$N_{G}(v)=\{x\in V(G):vx\in E(G)\}$. If $U$ is a nonempty subset of $V(G)$ 
then $\langle U \rangle$ denotes the subgraph of $G$ induced by $U$. 

 A graph $G$ is \textit{
hamiltonian }if it has a \textit{hamiltonian cycle }(a cycle containing all
the vertices of $G$), and \textit{traceable} if it has a \textit{hamiltonian
path }(a path containing all the vertices of $G$). A graph$\ G$ is \textit{
maximal nonhamiltonian}\emph{\ }(MNH)\ if $G$ is not hamiltonian, but $\ G+e$
is hamiltonian for each $e$ $\in E(\overline{G})$, where $\overline{G}$
denotes the complement of $G.$ A graph $G$ is \textit{maximal nontraceable}
(MNT) if $G$ is not traceable, but\ $G+e$ is traceable for each $e\in E(%
\overline{G}).$ A graph $G$ is \textit{hypohamiltonian}\ if $G$ is not
hamiltonian, but every vertex-deleted subgraph $G-v$ of $G$ is hamiltonian.
We say that a graph $G$ is \textit{maximal hypohamiltonian} (MHH) if it is
MNH and hypohamiltonian.

\bigskip

In 1978 Bollob\'{a}s \cite{bollobas} posed the problem of finding the least
number of edges, $f(n)$, in a MNH graph of order $n$. Bondy \cite{bondy} had
already shown that a MNH graph with order $n\geq 7$ that contained $m$
vertices of degree $2$ had at least $(3n+m)/2$ edges, and hence $f(n)\geq
\left\lceil 3n/2\right\rceil $ for $n\geq 7$. Combined results of Clark, Entringer and Shapiro 
\cite{ce}, \cite{ces} and Lin, Jiang, Zhang and Yang \cite
{ljzy} show that $f(n)=\left\lceil 3n/2\right\rceil $ for $n\geq 19$ and
for $n=6,10,11,12,13,17$. The values of $f(n)$ for the remaining values of $n$ are also given in \cite
{ljzy}.

\bigskip

Let $g(n)$ be the minimum size of a MNT graph of order $n$. Dudek, Katona
and Wojda \cite{dudek} showed that $g(n)\geq (3n-20)/2$ for all $n$ and,
by means of a recursive construction, they found MNT graphs of order $n$ and size 
$O(n\log n)$. To date, no cubic MNT graphs have been reported. We construct an
infinite family of cubic MNT graphs, thus showing that $g(n)\leq 3n/2$ for
infinitely many $n$.

Now let $g_{2}(n)$ be the minimum size of a 2-connected MNT graph of order $
n$. We prove that $g_{2}(n)\geq \left\lceil 3n/2\right\rceil $ for $n\geq 7$
. It then follows from our constructions that $g_{2}(n)=\left\lceil
3n/2\right\rceil $ for $n=8p$ for $p\geq 5$, $n=8p+2$ for $p\geq 6$, $n=8p+4$
for $p=3$ and $p \ge 6$, and $n=8p+6$ for $p\geq 4$.

\section{A lower bound for the size of a 2-connected MNT graph}

Bondy \cite{bondy} proved that if $G$ is a 2-connected MNH graph and $v\in
V(G)$ with degree $d\left( v\right) =2$, then each neighbour of $v$ has
degree at least 4. He also showed that the neighbours of such a vertex are
in fact adjacent.

\bigskip

In order to prove a corresponding result for 2-connected MNT graphs we need
the following result.

\begin{lemma}
\label{subgraph}Let $Q$ be a path in a MNT graph $G$. If $\langle V(Q) \rangle$ 
is not complete, then some internal vertex of $Q$ has a neighbour in $G-V(Q)$. 
\end{lemma}

\begin{proof}
Let $u$ and $v$ be two nonadjacent vertices of  $\langle V(Q) \rangle$. Then $G+uv$ has a hamiltonian 
path $P$. Let $x$ and $y$ be the two endvertices of $Q$ and suppose no internal vertex 
of $Q$ has a neighbour in $G-V(Q)$. Then $P$ has a subpath $R$ in 
$\langle V(Q) \rangle + uv$ and $R$ has either one or both endvertices in $\{x,y\}$. If 
$R$ has only one endvertex in $\{x,y\}$, then $P$ has an endvertex in $Q$. In either 
case the path obtained from $P$ by replacing $R$ with $Q$ is a hamiltonian path of $G$. 
\hfill \end{proof}

\begin{lemma}
\label{mntdeg2}If $G$ is a MNT graph and $v\in V(G)$ with $d\left( v\right)
=2$, then the neighbours of $v$ are adjacent. If in addition $G$ is
$2$-connected, then each neighbour of $v$ has degree at least $4$.
\end{lemma}

\begin{proof}
Let $N_G(v)=\{x_1,x_2\}$ and let $Q$ be the path $x_1vx_2$. Since $N_G(v)\subseteq Q$, 
it follows from Lemma \ref{subgraph} that $\langle V(Q) \rangle$ is a complete graph;
 hence $x_1$ and $x_2$ are adjacent.

Now assume that $G$ is $2$-connected. Since $G$ is not traceable we assume $
d(x_{1})>2$. Then also $d(x_{2})>2$ otherwise $x_{1}$ would be a cut vertex
of $G$.

Let $z$ be a neighbour of $x_1$ and let $Q$ be the path $zx_1vx_2$. Since $d(v)=2$ 
the graph $\langle V(Q) \rangle$ is not complete, and hence it follows from 
Lemma \ref{subgraph} that $x_1$ has a neighbour in $G-V(Q)$. Thus $d(x_1)\geq 4$. 
Similarly $d(x_2)\geq 4$.
\hfill \end{proof}

\bigskip

We also have the following two lemmas concerning MNT graphs that have
vertices of degree 2.

\begin{lemma}
\label{mnt1deg2}Suppose $G$ is a $2$-connected MNT graph. Suppose $v_{1},v_{2}\in
V(G)$  such that $d(v_{1})=d(v_{2})=2$ and $v_{1}$ and $
v_{2}$ have exactly one common neighbour $x$. Then $d(x)\geq 5.$
\end{lemma}

\begin{proof}
The vertices $v_{1}$ and $v_{2}$ cannot be adjacent otherwise $x$ would be a
cut vertex. Let $N(v_i)=\{x,y_i\}$; $i=1,2$. It follows from Lemma \ref{mntdeg2} that 
$x$ is adjacent to $y_i$; $i=1,2$. Let $Q$ be the path $y_1v_1xv_2y_2$. Since 
$\langle V(Q) \rangle$ is not complete, it follows from Lemma \ref{subgraph} that $x$ has 
a neighbour in $G-V(Q)$. Hence $d(x) \geq 5$.

\hfill \end{proof}

\begin{lemma}
\label{mnt2deg2}Suppose $G$ is a MNT graph. Suppose $v_{1},v_{2}\in V(G)$ such that
 $d(v_{1})=d(v_{2})=2$ and $v_{1}$ and $v_{2}$ have
the same two neighbours $x_{1}$ and $x_{2}$. Then $N_{G}(x_{1})-\{x_{2}
\}=N_{G}(x_{2})-\{x_{1}\}$. Also $d(x_{1})=d(x_{2})\geq 5.$
\end{lemma}

\begin{proof}
From Lemma \ref{mntdeg2} it follows that $x_{1}$ and $x_{2}$ are adjacent.
Let $Q$ be the path $x_2v_1x_1v_2$. $\langle V(Q) \rangle$ is not complete since $v_1$ 
and $v_2$ are not adjacent. Thus it follows from Lemma \ref{subgraph} that $x_1$ has a neighbour 
in $G-V(Q)$. Now suppose $p \in N_{G-V(Q)}(x_1)$ and  $p \notin N_{G}(x_2)$. 
Then a hamiltonian path $P$ in $G+px_{2}$ contains a
subpath of either of the forms given in the first column of Table 1. Note
that $i,j\in \{1,2\}$; $i\neq j$ and that $L$ represents a subpath of $P$ in $
G-\{x_{1},x_{2},v_{1},v_{2},p\}$. If each of the subpaths is replaced by the
corresponding subpath in the second column of the table we obtain a
hamiltonian path $P^{\prime }$ in $G$, which leads to a contradiction.

\begin{equation*}
\begin{tabular}{|l|l|}
\hline
Subpath of $P$ & Replace with \\ \hline
$v_{i}x_{1}v_{j}x_{2}p$ & $v_{i}x_{2}v_{j}x_{1}p$ \\ \hline
$v_{i}x_{1}Lpx_{2}v_{j}$ & $v_{i}x_{2}v_{j}x_{1}Lp$ \\ \hline
\end{tabular}
\end{equation*}
\begin{equation*}
\text{Table 1}
\end{equation*}
Hence $p \in N_{G}(x_2)$. Thus $N_{G}(x_{1})-\{x_{2}\}\subseteq N_{G}(x_{2})-\{x_{1}\}$. 
Similarly $N_{G}(x_{2})-\{x_{1}\}\subseteq N_{G}(x_{1})-\{x_{2}\}$. Thus $
N_{G}(x_{1})-\{x_{2}\}=N_{G}(x_{2})-\{x_{1}\}$ and hence $d(x_{1})=d(x_{2})$. Now let 
$Q$ be the path $px_1v_1x_2v_2$. Since $\langle V(Q) \rangle$ is not complete, 
it follows from Lemma \ref{subgraph} that $x_1$ or $x_2$ has a neighbour in $G-V(Q)$. 
Hence $d(x_1)=d(x_2) \geq 5$.
\hfill \end{proof}

\bigskip

We now consider the size of a 2-connected MNT\ graph.

\begin{lemma}
\label{mnt3deg2}Suppose $G$ is a MNT graph of order $n\geq 6$ and that $v_{1},v_{2}$ and $
v_{3}$ are vertices of degree $2$ in $G$ having the same neighbours, $x_{1}$
and $x_{2}$. Then $G-\{v_{1},v_{2},v_{3}\}$ is complete and hence $|E(G)|=
\frac{1}{2}(n^{2}-7n+24)$.
\end{lemma}

\begin{proof}
Suppose $G-\{v_{1},v_{2},v_{3}\}$ is not complete. Then there exist $p,q\in
V(G)-\{v_{1},v_{2},v_{3}\}$ which are not adjacent. However, since  $\{v_{1},v_{2}, 
v_{3}\}$ is an independent set, no path in $G+pq$ having $pq$
as an edge can contain all three of $v_{1},v_{2}$ and $v_{3}$. Thus $G+pq$
is not traceable. Thus $G-\{v_{1},v_{2},v_{3}\}=K_{n-3}.$ Hence $|E(G)|=
\frac{1}{2}(n-3)(n-4)+6.$
\hfill \end{proof}

\begin{theorem}
Suppose $G$ is a $2$-connected MNT graph. If $G$ has order $n\geq 7$ and $m$
vertices of degree $2$, then $|E(G)|\geq \frac{1}{2}(3n+m)$.
\end{theorem}

\begin{proof}
If $G$ has three vertices of degree 2 having the same two neighbours then, by 
Lemma \ref{mnt3deg2}
\begin{equation*}
|E(G)|=\tfrac{1}{2}(n^{2}-7n+24)\geq \tfrac{1}{2}(3n+m)\ \text{when}\ n \geq 7,
\end{equation*}
since $m=3.$

We now assume that $G$ does not have three vertices of degree 2 that have
the same two neighbours. Let $v_{1},...,v_{m}$ be the vertices of degree 2
in $G$ and let $H=G-\{v_{1},...,v_{m}\}.$ Then by Lemmas \ref{mntdeg2}, \ref
{mnt1deg2} and \ref{mnt2deg2} the minimum degree, $\delta (H)$ of $H$ is at
least 3. Hence 
\begin{equation*}
|E(G)|=|E(H)|+2m\geq \tfrac{3}{2}(n-m)+2m=\tfrac{1}{2}(3n+m).
\end{equation*}
\hfill \end{proof}

\bigskip

Thus $g_2(n) \geq \frac{1}{2} (3n+m)$ for $n \geq 7$. 
For $m \geq 1$ this bound is realized for $n=7$ (a Zelinka Type I
graph \cite{zelinka}) and $n=18$ (a graph constructed in \cite{bullock}). These graphs are depicted in Fig.\ 1. 
\[
\includegraphics[scale=0.5]{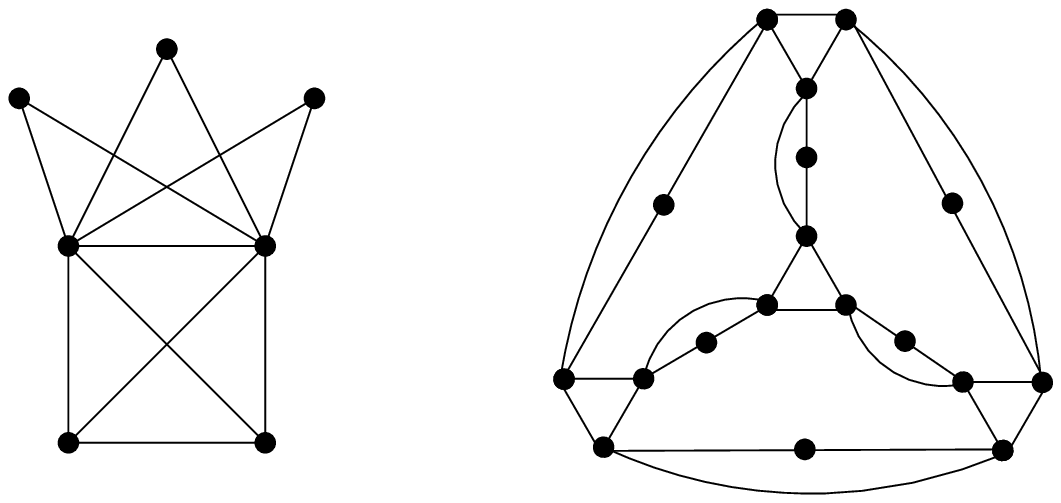}
\]
\[
\text{Fig.\ 1}
\]

We now construct an infinite family of 2-connected cubic MNT graphs, 
showing that $g_2(n)= \frac{3n}{2}$ for infinitely many $n$. 

\section{A construction of an infinite family of cubic MNT graphs}

In order to construct cubic MNT\ graphs we require the following lemmas
concerning MHH graphs.
\newpage
\begin{lemma}
\label{hypo2}Suppose $H$ is a hypohamiltonian graph having a vertex $z$ of
degree $3$. Put $F=H-z$.
\begin{enumerate}
\item[(a)] $F$ has a hamiltonian path ending at any of its vertices.
\item[(b)] There is no hamiltonian path in $F$ with both endvertices in $N_{H}(z)$.
\item[(c)] For any $y\in N_{H}(z)$ there exists a hamiltonian path in $F-y$ 
with the other two vertices of $N_{H}(z)$ being the endvertices.
\end{enumerate}
\end{lemma}

\begin{proof}
\begin{enumerate}
\item[(a)] $F$ is hamiltonian.
\item[(b)] If a hamiltonian path exists in $F$ having both endvertices in $N_{H}(z)$,
then $H$ has a hamiltonian cycle, which is a contradiction.
\item[(c)] Since $H-y$ is hamiltonian there is a hamiltonian cycle in $H-y$ containing 
the path $vzw$, where $v,w\in N_{H}(z)-\{y\}$. Thus there
is a hamiltonian path in $F-y$ with endvertices $v$ and $w$.
\end{enumerate}
\hfill \end{proof}

\begin{lemma}
\label{mnt1}Suppose $H$ is a MNH graph having a vertex $z$ of degree $3$. Put $F=H-z$.
If $u_{1}$ and $u_{2}$ are nonadjacent vertices in $F$, then $
F+u_{1}u_{2}$ has a hamiltonian path with both endvertices in $N_{H}(z)$.
\end{lemma}

\begin{proof}
There exists a hamiltonian cycle in $H+u_{1}u_{2}$ which contains the
path $vzw$, where $v,w\in N_{H}(z)$. Thus there exists a hamiltonian path
in $F+u_{1}u_{2}$ with endvertices $v$ and $w$.
\hfill \end{proof}

\bigskip

\begin{center}
\Large{\bf Construction of the graph $K_4[H_1,H_2,H_3]$}
\end{center}

For $i=1,2,3$, let $H_{i}$ be a cubic MHH graph, with a  vertex $z_i$ with neighbours $a_i$, $ b_i$ and $c_i$, which satisfies the following condition.\newline

\noindent{\textbf{Condition (C)}:}\newline
{\it{For every vertex}} $u_{i} \notin N_{H_{i}}(z_{i})$, {\it{the graph}} $H_{i}+z_{i}u_{i}$ 
{\it{has a hamiltonian cycle containing the edge}} $a_{i}z_{i}$ {\it{as well as a
hamiltonian cycle not containing}} $a_{i}z_{i}$.

\bigskip 

Graphs satisfying this condition will be presented at the end of the paper.

\bigskip

In the same sense as Gr\"{u}nbaum \cite{grunbaum} we use $H_{i}\backslash
z_{i}$ to denote $H_{i}$ ``opened up'' at $z_{i}$ (see Fig.\ 2). 
\[
\includegraphics[scale=0.67]{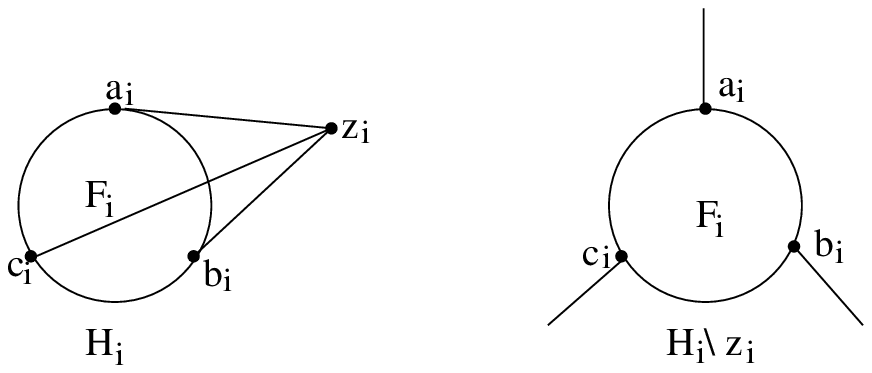}
\]
\[
\text{Fig.\ 2}
\]
\newpage
Let $K_4[H_1,H_2,H_3]$ be an inflated $K_4$ obtained from $H_{i}\backslash
z_{i}$; $i=1,2,3$ and a vertex $x$ by joining $x$ to the semi-edge incident 
with $a_i$ for $i=1,2,3$ and joining the remaining semi-edges as depicted in Fig.\ 3. Let $F_i$ denote 
$H_i-z_i$; $i=1,2,3$. We call $a_i$, $b_i$ and $c_i$ the \textit{exit vertices} of $F_i$.

\[
\includegraphics[scale=0.7]{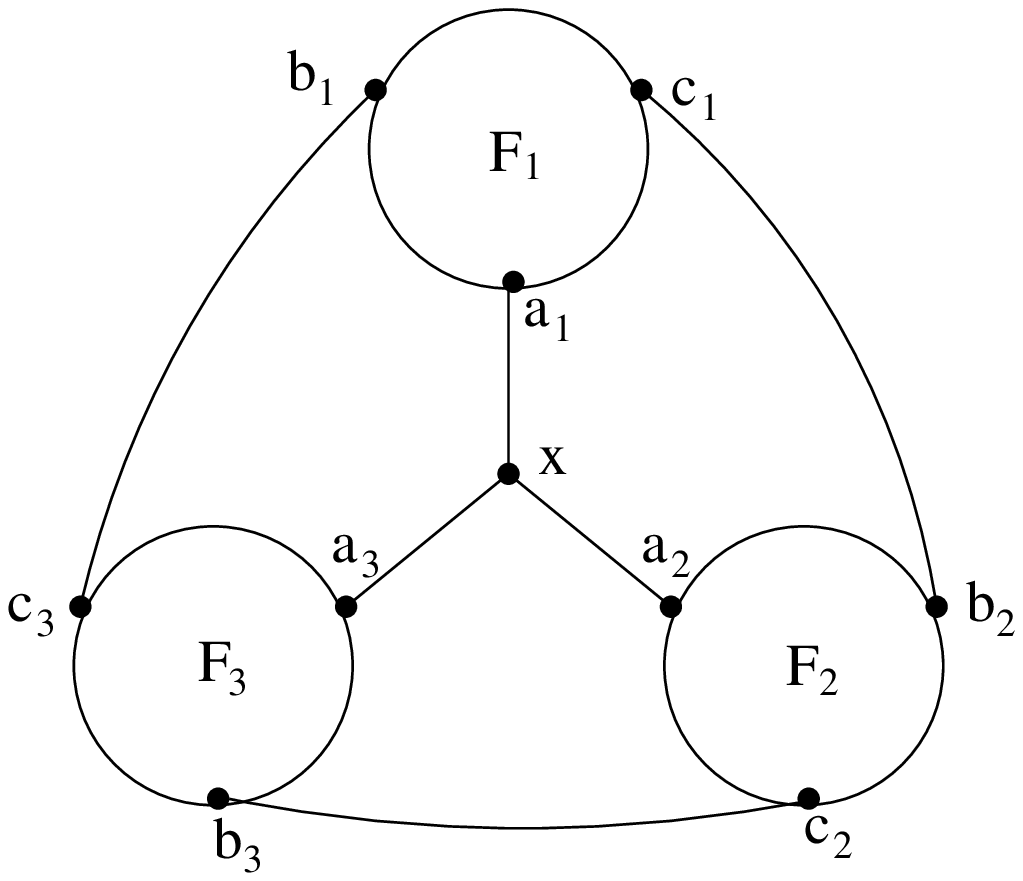}
\]
\[
\text{Fig.\ 3}
\]

We introduce the following notation which we use in the theorem below:
\newline
$P_{G}(v,w)$ denotes a hamiltonian path in a graph $G$ from $v$ to $w$;
\newline
$P_{G}(-,w)$ denotes a hamiltonian path in $G$ ending at $w$;\newline
$P_{G}(v,-)$ denotes a hamiltonian path in $G$ beginning at $v$; and\newline
$P_{G}$ denotes a hamiltonian path in $G$.

\begin{theorem}
\label{cubicmnt}The graph $G=K_{4}[H_{1},H_{2},H_{3}]$ is a cubic MNT graph.
\end{theorem}

\begin{proof}
It is obvious from the construction that $G$ is cubic.

We now show that $G$ is nontraceable. Suppose $P$ is a hamiltonian path of $
G $. Then at least one of the $F_{i}$'s, say $F_{2}$, does not contain an
endvertex of $P$. Thus $P$ passes through $F_{2}$, using two of the exit
vertices of $F_{2}$. However, by Lemma \ref{hypo2}(b) such a path cannot contain all the 
vertices of $F_{2}$.

We now show that $G+uv$ is traceable for all nonadjacent vertices $u$ and $v$
in $G$.\newline

\textbf{Case 1.} $u,v\in F_{i}$; $i\in \{1,2,3\}$.\newline
Without loss of generality consider $i=2$. By Lemma \ref{mnt1} there is
a hamiltonian path in $F_{2}+uv$ with endvertices two of $a_{2}$, $b_{2}$
and $c_{2}$.\newline
\textit{Subcase (i).} Suppose the endvertices are $a_{2}$ and $c_{2}$. (A
similar proof holds for $a_{2}$ and $b_{2}$.) By using Lemma \ref{hypo2}(a) we obtain the hamiltonian path 
\begin{equation*}
P_{G+uv}=P_{F_{1}}(-,a_{1})xP_{F_{2}+uv}(a_{2},c_{2})P_{F_{3}}(b_{3},-).
\end{equation*}
\textit{Subcase (ii).} Suppose the endvertices are $b_{2}$ and $c_{2}$. By
using  Lemma \ref{hypo2}(c) we obtain the hamiltonian path 
\begin{equation*}
P_{G+uv}=a_{1}xP_{F_{3}-c_{3}}(a_{3},b_{3})P_{F_{2}+uv}(c_{2},b_{2})P_{G_{1}-a_{1}}(c_{1},b_{1})c_{3}.
\end{equation*}
\textbf{Case 2.} $u\in \{a_{i},b_{i},c_{i}\}$ and $v\in
\{a_{j},b_{j},c_{j}\}$; $i,j\in \{1,2,3\}$; $i\neq j $.\newline
Without loss of generality we choose $i=2$ and $j=3$. By using  Lemmas \ref{hypo2}(a) and 
(c) we find a hamiltonian path $P_{G+uv}$ in $G+uv$. All
subcases can be reduced to the following:\newline
\textit{Subcase (i).} $u=a_{2},v=a_{3}$. 
\begin{equation*}
P_{G+uv}=a_{2}a_{3}xP_{G_{1}-c_{1}}(a_{1},b_{1})P_{F_{3}-a_{3}}(c_{3},b_{3})P_{F_{2}-a_{2}}(c_{2},b_{2})c_{1}.
\end{equation*}
\textit{Subcase (ii).} $u=a_{2},v=b_{3}$. 
\begin{equation*}
P_{G+uv}=c_{2}b_{3}P_{F_{2}-c_{2}}(a_{2},b_{2})P_{G_{1}-b_{1}}(c_{1},a_{1})xP_{F_{3}-b_{3}}(a_{3},c_{3})b_{1}.
\end{equation*}
\textit{Subcase (iii).} $u=a_{2},v=c_{3}$. 
\begin{equation*}
P_{G+uv}=b_{1}c_{3}P_{F_{2}-c_{2}}(a_{2},b_{2})P_{G_{1}-b_{1}}(c_{1},a_{1})xP_{F_{3}-c_{3}}(a_{3},b_{3})c_{2}.
\end{equation*}
\textit{Subcase (iv).} $u=b_{2},v=b_{3}$. 
\begin{equation*}
P_{G+uv}=c_{2}b_{3}P_{F_{2}-c_{2}}(b_{2},a_{2})xP_{F_{3}-b_{3}}(a_{3},c_{3})P_{F_{1}}(b_{1},-).
\end{equation*}
\textit{Subcase (v).} $u=b_{2},v=c_{3}$. 
\begin{equation*}
P_{G+uv}=c_{1}b_{2}c_{3}P_{G_{1}-c_{1}}(b_{1},a_{1})xP_{F_{2}-b_{2}}(a_{2},c_{2})P_{F_{3}-c_{3}}(b_{3},a_{3}).
\end{equation*}
\textbf{Case 3.} $u\in F_{i}-\{a_{i},b_{i},c_{i}\}$ and $v\in F_{j}$; $i,j\in \{1,2,3\}$; 
$i\neq j$.\newline
Without loss of generality we choose $i=2$ and $j=3$. Let $F_{2}^{\ast}$ be the graph 
obtained from $G$ by contracting $G-V(F_2)$ to a single vertex $z_{2}^{\ast}$. Then $F_{2}^{\ast}$ 
is isomorphic to $H_2$ and hence, by Condition  \textbf{(C)}, $F_{2}^{\ast}+uz_{2}^{\ast}$ has a hamiltonian 
cycle containg the path $uz_{2}^{\ast}a_2$. Thus $F_2$ has 
a hamiltonian path  with endvertices $u$ and $a_{2}$. 
Using this fact and Lemma \ref{hypo2}(a) we construct the hamiltonian path 
\begin{equation*}
P_{G+uv}=P_{F_{3}}(-,v)P_{F_{2}}(u,a_{2})xP_{F_{1}}(a_{1},-).
\end{equation*}
\textbf{Case 4.} $u=x$ and $v\in F_{i}$; $i\in \{1,2,3\}$.\newline
Without loss of generality we choose $i=2$.\newline
\textit{Subcase (i).} $v\in \{b_{2},c_{2}\}$.\newline
Consider $v=b_{2}.$ (The case $v=c_{2}$ follows similarly.) By using  Lemmas \ref{hypo2}(a) and 
(c) we obtain the hamiltonian path 
\begin{equation*}
P_{G+uv}=P_{F_{3}}(-,b_{3})P_{F_{2}-b_{2}}(c_{2},a_{2})xb_{2}P_{F_{1}}(c_{1},-).
\end{equation*}
\textit{Subcase (ii).} $v\in F_{2}-\{a_{2},b_{2},c_{2}\}$.\newline
According to Condition  \textbf{(C)} and an argument similar to that in Case 3, there is 
a hamiltonian path in $F_{2}$ with
endvertices $v$ and $d$, where $d\in \{b_{2},c_{2}\}$. Suppose $d=b_{2}$. (A
similar proof holds for $d=c_{2}$.) Using this fact and Lemma \ref{hypo2}(a) we
construct the hamiltonian path  
\begin{equation*}
P_{G+uv}=P_{F_{3}}(-,a_{3})xP_{F_{2}}(v,b_{2})P_{F_{1}}(c_{1},-).
\end{equation*}
\hfill \end{proof}

\bigskip

The Petersen graph ($n=10$), the Coxeter graph ($n=28$) and the Isaacs'
snarks $J_{k}$ ($n=4k$) for odd $k\geq 5$ are all cubic MHH graphs (see \cite
{bondy}, \cite{ce}). We determined, by using the Graph Manipulation Package developed by Siqinfu 
and Sheng Bau*, that a snark of order $22$, reported by Chisala \cite{chisala}, is also MHH.
The Petersen graph, the snark of order 22 and the Coxeter graph are shown in Fig.\ 4.
\[
\includegraphics[scale=0.6]{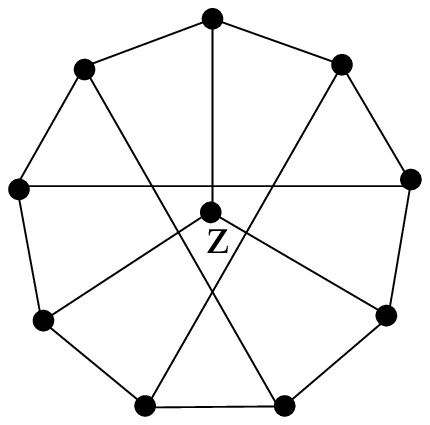}
\ \ \ \ \ 
\includegraphics[scale=0.5]{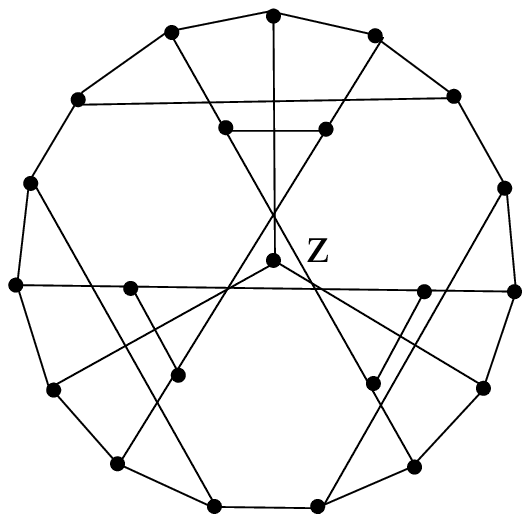}
\ \ \ \ \
\includegraphics[scale=0.4]{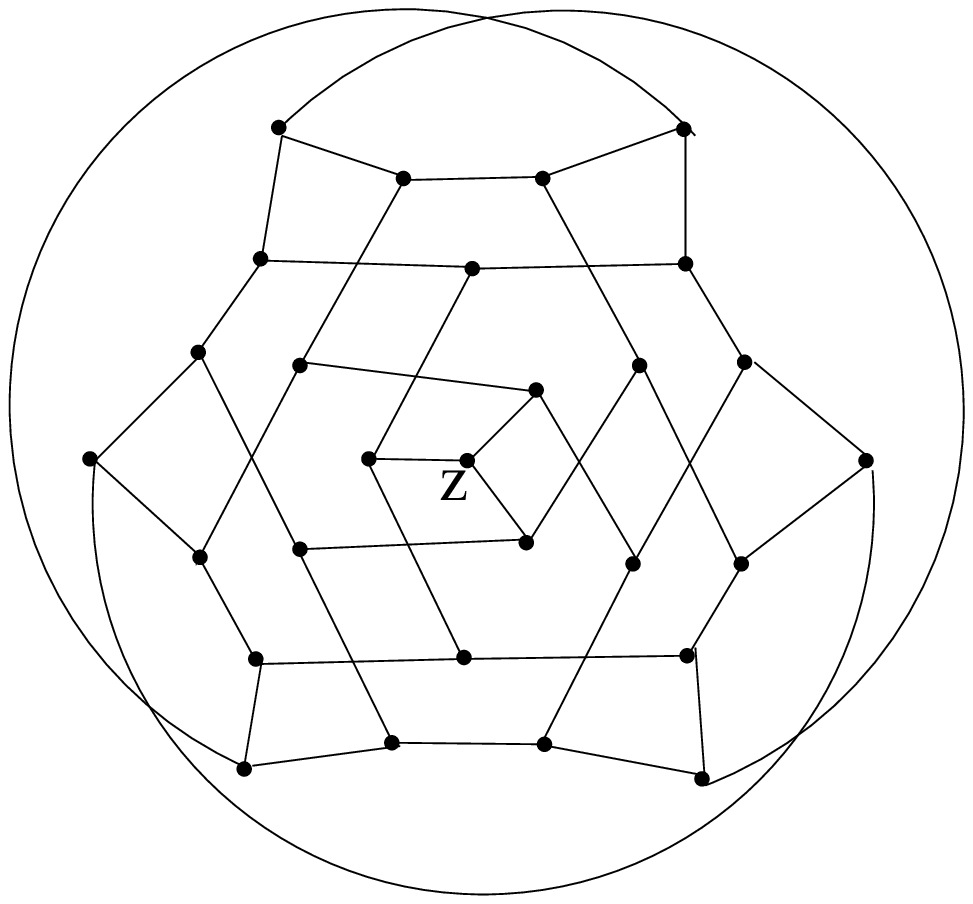}
\]
\[
\text{Fig.\ 4}
\]

All the cubic MHH graphs mentioned above satisfy condition \textbf{(C)}. In fact, it
follows from Theorem 10 in \cite{ce} that the Isaac's snarks satisfy a
stronger condition, namely that if $v\notin N_H(z)$ then, for \textit{every} $u\in N_H(z)$ 
there exists a hamiltonian cycle in $H+uv$ containing $uz$, and this condition holds for 
all $z\in H$. We determined, again by using the Graph Manipulation Package, that each of the  
graphs shown in Fig.\ 4 also satisfy this extended condition for the specified vertex $z$. 
The programme allows one to sketch a graph $G$ on the computer screen by placing vertices and 
adding edges. On request the programme will either draw in a hamiltonian cycle or state that the graph 
is non-hamiltonian. We tested to see if each of the above mentioned MNH graphs $G$ satisfied the 
condition \textbf{(C)} by considering symmetry and adding appropriate edges and 
noting the structure of the hamiltonian cycle drawn in.

Thus, by using various combinations of these MHH graphs, we can
produce cubic MNT\ graphs of order 
\begin{equation*}
n=\left\{ 
\begin{tabular}{ll}
$8p$ & $p\geq 5$ \\ 
$8p+2$ & $p\geq 6$ \\ 
$8p+4$ & $p=3,p\geq 6$** \\ 
$8p+6$ & $p\geq 4.$
\end{tabular}
\right.
\end{equation*}

Thus $g_{2}(n)= \frac{3n}{2} $ for all the values of $
n$ stated above.

\bigskip

\textbf{Remark: }Our construction yields MNT graphs of girths 5, 6 and 7. We
do not know whether MNT graphs of girth bigger than 7 exist.

\bigskip

\textbf{Acknowledgements} \newline
* We wish to thank Sheng Bau for allowing us the use of the 
programme, Graph Manipulation Package  Version 1.0 (1996), Siqinfu and Sheng Bau,
Inner Mongolia Institute of Finance and Economics, Huhhot, CN-010051, People's 
Republic of China. \newline
** We wish to thank the referee who brought to our attention the infinite family 
$K_{4}[S,J_{k},J_{k'}]$ of MNT graphs, where $S$ is the snark of order $22$ and $J_{k}$ and $J_{k'}$ 
are Isaac's snarks, which gives $n=8p+4$ for $p \ge7$.

\end{document}